\documentclass{article}
\usepackage{amssymb}
\usepackage{amsfonts}
\usepackage{amsmath}

\newtheorem{theorem}{Theorem}

\newtheorem{corollary}[theorem]{Corollary}

\newtheorem{definition}[theorem]{Definition}

\input{tcilatex}
\begin{document}

\title{Intrinsic Equations For a Relaxed Elastic Line of Second Kind on an
Oriented Surface}
\author{Ergin Bayram and Emin Kasap \\
Ondokuz Mayis University, Faculty of Arts and Sciences, \\
Mathematics Department, Samsun, Turkey}
\maketitle

\begin{abstract}
Let $\alpha \left( s\right) $ be an arc on a connected oriented surface $S$
in $E^{3}$, parameterized by arc length $s$, with torsion $\tau \ $and
length $l$. The total square torsion $F$ of $\alpha $ is defined by $%
F=\int_{0}^{l}\tau ^{2}ds$. The arc $\alpha $ is called a relaxed elastic
line of second kind if it is an extremal for the variational problem of
minimizing the value of $F$ within the family of all arcs of length $l$ on $S
$ having the same initial point and initial direction as $\alpha $. In this
study, we obtain differential equation and boundary conditions for a relaxed
elastic line of second kind on an oriented surface.
\end{abstract}

\section{Introduction}

Let $\alpha \left( s\right) ~$denote an arc on a connected oriented surface $%
S$ in $E^{3},~$parameterized by arc length $s$, $0\leq s\leq l,$ with
curvature $\kappa \left( s\right) .$ The total square curvature $K$ of $%
\alpha $ is defined by

\bigskip 
\begin{equation}
K=\int_{0}^{l}\kappa ^{2}ds.  \label{1.1}
\end{equation}

An arc is called a relaxed elastic line if it is an extremal for the
variational problem of minimizing the value of $K~$within the family of all
arcs of length $l$ on $S$ having the same initial point and initial
direction as $\alpha $ \cite{manning}. In \cite{manning} they derive the
intrinsic equations for a relaxed elastic line on an oriented surface.
Hilbert and Cohn-Vossen \cite{hilbert} incorrectly suggested a flexible
knitting needle, constrained to conform to a surface, as one model for a
geodesic on a surface. This model actually gives a relaxed elastic line on
the surface, and is not generally a geodesic unless the surface lies in a
plane or on a sphere \cite{manning}. Physical motivation for study of the
problem of elastic lines\ on surfaces may be found in the nucleosome core
particle \cite{nucleosome 1}, \cite{nucleosome 2}.

There are several papers about this kind of minimization problems \cite%
{yucesan}, \cite{ekici}, \cite{yucesan2}, \cite{sahin}.

In \cite{zuhal} authors handled the problem of minimizing the total square
torsion on an oriented surface and defined the relaxed elastic line of
second kind. However, they only gave Euler-Lagrange equations for this
problem.

In this paper, we obtain intrinsic equations for a relaxed elastic line of
second kind. We give a differential equation and three boundary conditions.

\section{Derivation of equations}

Let $\alpha \left( s\right) ~$denote an arc on a connected oriented surface $%
S$ in $E^{3},~$parameterized by arc length $s$, $0\leq s\leq l,$ with
torsion $\tau \left( s\right) .$ The total square torsion $F$ of $\alpha $
is defined by

\bigskip 
\begin{equation}
F=\int_{0}^{l}\tau ^{2}ds.  \label{2.1}
\end{equation}

\begin{definition}
The arc $\alpha $ is called a relaxed elastic line of second kind if it is
an extremal for the variational problem of minimizing the value of $F~$%
within the family of all arcs of length $l$ on $S$ having the same initial
point and initial direction as $\alpha \ \cite{zuhal}.$
\end{definition}

We assume that the arc $\alpha $ is smooth enough to have derivatives up to
required order and parameterized with arc length. Also, for technical
reasons we assume that $\kappa \neq 0,\ \forall s\ $on $\alpha $. On $\alpha 
$, let $T\left( s\right) =\alpha ^{\prime }\left( s\right) $ denote the unit
tangent vector field$,$ $n\left( s\right) $ denote the unit surface normal
vector field to $S\ $and $Q\left( s\right) =nxT.\ $Then, $\left\{
T,Q,n\right\} $ gives an orthonormal basis on $\alpha $ and $\left\{
T,Q\right\} $ gives a basis for the vectors tangent to $S$ at $\alpha \left(
s\right) .$ The frame $\left\{ T,Q,n\right\} $ is called Darboux frame.
Derivative equations for the Darboux frame is

\begin{equation}
\left( 
\begin{array}{c}
T^{\prime } \\ 
Q^{\prime } \\ 
n^{\prime }%
\end{array}%
\right) =\left( 
\begin{array}{c}
0\ \ \ \ \ \ \ \kappa _{g}\ ~\ \ \kappa _{n} \\ 
-\kappa _{g}\ \ \ \ \ \ 0~\ \ \ \ \ \tau _{g} \\ 
-\kappa _{n}~-\tau _{g}\ \ \ \ ~0%
\end{array}%
\right) \left( 
\begin{array}{c}
T \\ 
Q \\ 
n%
\end{array}%
\right) ,  \label{2.2}
\end{equation}%
where $\kappa _{g},\ \kappa _{n}\ $and $\tau _{g}\ $are geodesic curvature,
normal curvature and geodesic torsion, respectively \cite{carmo}. The square
curvature $\kappa ^{2}$ and the torsion $\tau $ of $\alpha $ on S are given
by%
\begin{equation}
\left\{ 
\begin{array}{c}
\kappa ^{2}=\kappa _{g}^{2}+\kappa _{n}^{2}, \\ 
\tau =\tau _{g}+\frac{\kappa _{g}\kappa _{n}^{\prime }-\kappa _{g}^{\prime
}\kappa _{n}}{\kappa _{g}^{2}+\kappa _{n}^{2}}.%
\end{array}%
\right.  \label{2.3}
\end{equation}

Suppose that $\alpha $ lies in a coordinate patch $\left( u,v\right)
\rightarrow x\left( u,v\right) \ $of $S,\ $and let $x_{u}=\partial
x/\partial u,\ x_{v}=\partial x/\partial v.\ $Then, $\alpha $ is expressed as

\begin{equation*}
\alpha \left( s\right) =x\left( u\left( s\right) ,v\left( s\right) \right)
,\ 0\leq s\leq l,
\end{equation*}%
with

\begin{equation*}
T\left( s\right) =\alpha ^{\prime }\left( s\right) =\frac{du}{ds}x_{u}+\frac{%
dv}{ds}x_{v}
\end{equation*}%
and

\begin{equation*}
Q\left( s\right) =p\left( s\right) x_{u}+q\left( s\right) x_{v}
\end{equation*}%
for suitable scalar functions $p\left( s\right) \ $and $q\left( s\right) .$

Now, we will define variational fields for our problem. In order to obtain
variational arcs of length $l$, we need to extend $\alpha $ to an arc $%
\alpha ^{\ast }\left( s\right) $ defined for $0\leq s\leq l^{\ast },$ with $%
l^{\ast }>l$ but sufficiently close to $l$ so that $\alpha ^{\ast }$ lies in
the coordinate patch. Let $\mu \left( s\right) ,\ 0\leq s\leq l^{\ast },\ $%
be a scalar function of class $C^{2},\ $not vanishing identically. Define%
\begin{equation*}
\eta \left( s\right) =\mu \left( s\right) p^{\ast }\left( s\right) ,\ \
\zeta \left( s\right) =\mu \left( s\right) q^{\ast }\left( s\right) .
\end{equation*}

Then,

\begin{equation}
\eta \left( s\right) x_{u}+\zeta \left( s\right) x_{v}=\mu \left( s\right)
Q\left( s\right)  \label{2.4}
\end{equation}%
along $\alpha .\ $Also assume that

\begin{equation}
\mu \left( 0\right) =0,\ \mu ^{\prime }\left( 0\right) =0,~\mu ^{\prime
\prime }\left( 0\right) =0.  \label{2.5}
\end{equation}

Now, define

\begin{equation}
\beta \left( \sigma ;t\right) =x\left( u\left( \sigma \right) +t\eta \left(
\sigma \right) ,v\left( \sigma \right) +t\zeta \left( \sigma \right) \right)
,  \label{2.6}
\end{equation}%
for $0\leq \sigma \leq $\bigskip $l^{\ast }.$ For $\left\vert t\right\vert
<\varepsilon _{1}\ $(where $\varepsilon _{1}>0\ $depends upon the choice of $%
\alpha ^{\ast }$ and of $\mu $), the point $\beta \left( \sigma ;t\right) $
lies in the coordinate patch. For fixed $t,\ \beta \left( \sigma ;t\right) \ 
$gives an arc with the same initial point and initial direction as $\alpha
,\ $because of $\left( \ref{2.5}\right) $. For $t=0,$ $\beta \left( \sigma
;0\right) \ $is the same as $\alpha ^{\ast }\ $and $\sigma \ $is arc length.
For $t\neq 0,\ $the parameter $\sigma \ $is not arc length in general.

For fixed $t,\ \left\vert t\right\vert <\varepsilon _{1},$ let $L^{\ast
}\left( t\right) $ denote the length of the arc $\beta \left( \sigma
;t\right) $, $0\leq \sigma \leq l^{\ast }.$ Then,

\begin{equation}
L^{\ast }\left( t\right) =\int_{0}^{l}\sqrt{\left\langle \frac{\partial
\beta }{\partial \sigma },\frac{\partial \beta }{\partial \sigma }%
\right\rangle }d\sigma  \label{2.7}
\end{equation}%
with

\begin{equation}
L^{\ast }\left( 0\right) =l^{\ast }>l.  \label{2.8}
\end{equation}

\bigskip By $\left( \ref{2.6}\right) \ $and $\left( \ref{2.7}\right) \
L^{\ast }\left( t\right) \ $is continuous and differentiable in $t.\ $%
Particularly, it follows from $\left( \ref{2.8}\right) $ that

\begin{equation}
L^{\ast }\left( t\right) >\frac{l+l^{\ast }}{2}>l\ \ \ \text{for \ }%
\left\vert t\right\vert <\varepsilon  \label{2.9}
\end{equation}%
for a suitable $\varepsilon \ $satisfying $0<\varepsilon \leq \varepsilon
_{1}.\ $Because of $\left( \ref{2.9}\right) \ $one can restrict $\beta
\left( \sigma ;t\right) ,\ 0\leq \left\vert t\right\vert <\varepsilon ,\ $to
an arc of length $l$ by restricting the parameter $\sigma $ to an interval $%
0\leq \sigma \leq \lambda \left( t\right) \leq l^{\ast }$ by requiring

\begin{equation}
\int_{0}^{\lambda \left( t\right) }\sqrt{\left\langle \frac{\partial \beta }{%
\partial \sigma },\frac{\partial \beta }{\partial \sigma }\right\rangle }%
d\sigma =l.  \label{2.10}
\end{equation}

Note that $\lambda \left( 0\right) =l.$ The function $\lambda \left(
t\right) $ need not be determined explicitly, but we shall need

\begin{equation}
\left. \frac{d\lambda }{dt}\right\vert _{t=0}=\int_{0}^{l}\mu \kappa _{g}ds.
\label{2.11}
\end{equation}

The proof of $\left( \ref{2.11}\right) $ and of other results will depend on
calculations from $\left( \ref{2.6}\right) $ such as

\begin{equation}
\left. \frac{\partial \beta }{\partial \sigma }\right\vert _{t=0}=T,\ \
0\leq s\leq l,  \label{2.12}
\end{equation}%
which gives

\begin{equation}
\left. \frac{\partial ^{2}\beta }{\partial \sigma ^{2}}\right\vert
_{t=0}=T^{\prime }=\kappa _{g}Q+\kappa _{n}n.  \label{2.13}
\end{equation}

Also

\begin{equation}
\left. \frac{\partial \beta }{\partial t}\right\vert _{t=0}=\mu Q
\label{2.14}
\end{equation}%
because of $\left( \ref{2.4}\right) $. Further differentiation of $\left( %
\ref{2.14}\right) $ gives

\begin{equation}
\left. \frac{\partial ^{2}\beta }{\partial t\partial \sigma }\right\vert
_{t=0}=\left. \frac{\partial ^{2}\beta }{\partial \sigma \partial t}%
\right\vert _{t=0}=\mu ^{\prime }Q+\mu Q^{\prime }=-\mu \kappa _{g}T+\mu
^{\prime }Q+\mu \tau _{g}n  \label{2.15}
\end{equation}%
and using $\left( \ref{2.2}\right) $,

\begin{eqnarray}
\left. \frac{\partial ^{3}\beta }{\partial t\partial \sigma ^{2}}\right\vert
_{t=0} &=&\left( -2\mu ^{\prime }\kappa _{g}-\mu \kappa _{g}^{\prime }-\mu
\kappa _{n}\tau _{g}\right) T+\left( \mu ^{\prime \prime }-\mu \kappa
_{g}^{2}-\mu \tau _{g}^{2}\right) Q  \label{2.16} \\
&&+\left( 2\mu ^{\prime }\tau _{g}+\mu \tau _{g}^{\prime }-\mu \kappa
_{g}\kappa _{n}\right) n.  \notag
\end{eqnarray}

Also using $\left( \ref{2.13}\right) $ we have

\begin{equation}
\left. \frac{\partial ^{3}\beta }{\partial \sigma ^{3}}\right\vert
_{t=0}=-\left( \kappa _{g}^{2}+\kappa _{n}^{2}\right) T+\left( \kappa
_{g}^{\prime }+\kappa _{n}\tau _{g}\right) Q+\left( \kappa _{n}^{\prime
}+\kappa _{g}\tau _{g}\right) n  \label{2.17}
\end{equation}%
and by $\left( \ref{2.14}\right) $

\begin{eqnarray}
\left. \frac{\partial ^{4}\beta }{\partial t\partial \sigma ^{3}}\right\vert
_{t=0} &=&\left( -3\mu ^{\prime \prime }\kappa _{g}+3\mu ^{\prime }\kappa
_{g}^{\prime }+\mu \kappa _{g}^{\prime \prime }+\mu \kappa _{n}^{\prime
}\tau _{g}+2\mu \kappa _{n}\tau _{g}^{\prime }\right.  \label{2.18} \\
&&\left. +3\mu ^{\prime }\kappa _{n}\tau _{g}-\mu \kappa _{g}\kappa
_{n}^{2}-\mu \kappa _{g}^{3}-\mu \kappa _{g}\tau _{g}^{2}\right) T  \notag \\
&&-\left( 3\mu ^{\prime }\kappa _{g}^{2}+3\mu ^{\prime }\tau _{g}^{2}+3\mu
\kappa _{g}\kappa _{g}^{\prime }+3\mu \tau _{g}\tau _{g}^{\prime }-\mu
^{\prime \prime \prime }\right) Q  \notag \\
&&-\left( 3\mu ^{\prime }\kappa _{g}\kappa _{n}+2\mu \kappa _{g}^{\prime
}\kappa _{n}+\mu \kappa _{g}^{2}\tau _{g}-3\mu ^{\prime \prime }\tau
_{g}\right.  \notag \\
&&\left. -3\mu ^{\prime }\tau _{g}^{\prime }-\mu \tau _{g}^{\prime \prime
}+\mu \kappa _{g}\kappa _{n}^{\prime }+\mu \tau _{g}^{3}+\mu \kappa
_{n}^{2}\tau _{g}\right) n.  \notag
\end{eqnarray}

Now, let $F\left( t\right) $ denote the functional of a relaxed elastic line
of second kind for the arc $\beta \left( \sigma ;t\right) ,\ 0\leq \sigma
\leq \lambda \left( t\right) ,\ \left\vert t\right\vert <\varepsilon .$
Since, in general, $\sigma \ $is not the arc length for $t\neq 0$ functional 
$\left( \ref{2.1}\right) $ can be calculated as follows:%
\begin{equation*}
F\left( t\right) =\int_{0}^{\lambda \left( t\right) }\left( \frac{%
\left\langle \frac{\partial \beta }{\partial \sigma }\times \frac{\partial
^{2}\beta }{\partial \sigma ^{2}},\frac{\partial ^{3}\beta }{\partial \sigma
^{3}}\right\rangle }{\left\langle \frac{\partial \beta }{\partial \sigma },%
\frac{\partial \beta }{\partial \sigma }\right\rangle \left\langle \frac{%
\partial ^{2}\beta }{\partial \sigma ^{2}},\frac{\partial ^{2}\beta }{%
\partial \sigma ^{2}}\right\rangle -\left\langle \frac{\partial ^{2}\beta }{%
\partial \sigma ^{2}},\frac{\partial \beta }{\partial \sigma }\right\rangle
^{2}}\right) ^{2}d\sigma .
\end{equation*}

A necessary condition for $\alpha $ to be an extremal is that

\begin{equation*}
\left. \frac{dF}{dt}\right\vert _{t=0}=0
\end{equation*}%
for arbitrary $\mu $ satisfying $\left( \ref{2.5}\right) $. In calculating $%
dF/dt$; we give explicitly only terms that do not vanish for $t=0$. The
omitted terms are those with the factor

\begin{equation*}
\left\langle \frac{\partial ^{2}\beta }{\partial \sigma ^{2}},\frac{\partial
\beta }{\partial \sigma }\right\rangle
\end{equation*}%
which vanishes at $t=0$ because $\left\langle T,T^{\prime }\right\rangle =0.$
Thus, using $\left( \ref{2.3}\right) $, $\left( \ref{2.11}-\ref{2.13}\right) 
$ and $\left( \ref{2.15}-\ref{2.18}\right) $ we get

\bigskip 
\begin{eqnarray*}
\left. \frac{dF}{dt}\right\vert _{t=0} &=&\int_{0}^{l}\mu \left\{ \kappa
_{g}\tau ^{2}\left( l\right) +2\frac{\tau }{\kappa ^{2}}\left[ -\kappa
_{g}\kappa ^{2}\tau +\kappa _{g}\tau _{g}\kappa ^{2}+\left( \kappa
_{g}\kappa _{n}-\tau _{g}^{\prime }\right) \left( \kappa _{g}^{\prime
}-\kappa _{n}\tau _{g}\right) \right. \right. \\
&&\left. -\left( \kappa _{g}^{2}+\tau _{g}^{2}\right) \left( \kappa
_{n}^{\prime }+\kappa _{g}\tau _{g}\right) -\kappa _{g}\left( 2\kappa
_{g}^{\prime }\kappa _{n}+\kappa _{g}^{2}\tau _{g}-\tau _{g}^{\prime \prime
}+\kappa _{g}\kappa _{n}^{\prime }+\tau _{g}^{3}+\kappa _{n}^{2}\tau
_{g}\right) \right. \\
&&\left. \left. +3\kappa _{n}\left( \kappa _{g}\kappa _{g}^{\prime }+\tau
_{g}\tau _{g}^{\prime }\right) -2\tau \left( -\kappa _{g}\kappa ^{2}+\kappa
_{n}\left( \tau _{g}^{\prime }-\kappa _{g}\kappa _{n}\right) -\kappa
_{g}\left( \kappa _{g}^{2}+\tau _{g}^{2}\right) \right) \right] \right\} ds
\\
&&+2\int_{0}^{l}\mu ^{\prime }\frac{\tau }{\kappa ^{2}}\left[ -\kappa
_{n}\kappa ^{2}-2\tau _{g}\left( \kappa _{g}^{\prime }-\kappa _{n}\tau
_{g}\right) -3\kappa _{g}\left( \kappa _{g}\kappa _{n}-\tau _{g}^{\prime
}\right) +3\kappa _{n}\left( \kappa _{g}^{2}+\tau _{g}^{2}\right) \right. \\
&&\left. -4\kappa _{n}\tau _{g}\tau \right] ds+2\int_{0}^{l}\mu ^{\prime
\prime }\frac{\tau }{\kappa ^{2}}\left( \kappa _{n}^{\prime }+4\kappa
_{g}\tau _{g}-\kappa _{g}\tau \right) ds-2\int_{0}^{l}\mu ^{\prime \prime
\prime }\frac{\kappa _{n}\tau }{\kappa ^{2}}ds.
\end{eqnarray*}

\bigskip However using integration by parts and $\left( \ref{2.15}\right) $
we have\bigskip 
\begin{eqnarray*}
\left. \frac{dF}{dt}\right\vert _{t=0} &=&\int_{0}^{l}\mu \left\{ \kappa
_{g}\tau ^{2}\left( l\right) +2\frac{\tau }{\kappa ^{2}}\left[ -\kappa
_{g}\kappa ^{2}\tau +\kappa _{g}\tau _{g}\kappa ^{2}+\left( \kappa
_{g}\kappa _{n}-\tau _{g}^{\prime }\right) \left( \kappa _{g}^{\prime
}-\kappa _{n}\tau _{g}\right) \right. \right. \\
&&\left. -\left( \kappa _{g}^{2}+\tau _{g}^{2}\right) \left( \kappa
_{n}^{\prime }+\kappa _{g}\tau _{g}\right) -\kappa _{g}\left( 2\kappa
_{n}\kappa _{g}^{\prime }+\kappa _{g}^{2}\tau _{g}-\tau _{g}^{\prime \prime
}+\kappa _{g}\kappa _{n}^{\prime }+\tau _{g}^{3}+\kappa _{n}^{2}\tau
_{g}\right) \right. \\
&&\left. +3\kappa _{n}\left( \kappa _{g}\kappa _{g}^{\prime }+\tau _{g}\tau
_{g}^{\prime }\right) +2\tau \left( \kappa _{g}\kappa ^{2}-\kappa _{n}\left(
\tau _{g}^{\prime }-\kappa _{g}\kappa _{n}\right) +\kappa _{g}\left( \kappa
_{g}^{2}+\tau _{g}^{2}\right) \right) \right] \\
&&-2\left[ \frac{\tau }{\kappa ^{2}}\left[ -\kappa ^{2}\kappa _{n}-2\tau
_{g}\left( \kappa _{g}^{\prime }-\kappa _{n}\tau _{g}\right) -3\kappa
_{g}\left( \kappa _{g}\kappa _{n}-\tau _{g}^{\prime }\right) +3\kappa
_{n}\left( \kappa _{g}^{2}+\tau _{g}^{2}\right) \right. \right. \\
&&\left. \left. \left. -4\kappa _{n}\tau _{g}\tau \right] \right] ^{\prime
}+2\left[ \frac{\tau }{\kappa ^{2}}\left[ \kappa _{n}^{\prime }+\kappa
_{n}\tau _{g}+3\kappa _{g}\tau _{g}-\kappa _{g}\tau \right] \right] ^{\prime
\prime }+2\left( \frac{\kappa _{n}\tau }{\kappa ^{2}}\right) ^{\prime \prime
\prime }\right\} ds \\
&&+\mu \left( l\right) \left\{ 2\frac{\tau \left( l\right) }{\kappa
^{2}\left( l\right) }\left[ -\kappa _{n}\left( l\right) \kappa ^{2}\left(
l\right) -2\tau _{g}\left( l\right) \left( \kappa _{g}^{\prime }\left(
l\right) -\kappa _{n}\left( l\right) \tau _{g}\left( l\right) \right)
\right. \right. \\
&&\left. -3\kappa _{g}\left( l\right) \left( \kappa _{g}\left( l\right)
\kappa _{n}\left( l\right) -\tau _{g}^{\prime }\left( l\right) \right)
+3\kappa _{n}\left( l\right) \left( \kappa _{g}^{2}\left( l\right) +\tau
_{g}^{2}\left( l\right) \right) -4\kappa _{n}\left( l\right) \tau _{g}\left(
l\right) \tau \left( l\right) \right] \\
&&\left. \left. -2\left( \frac{\tau }{\kappa ^{2}}\left[ \kappa _{n}^{\prime
}+3\kappa _{g}\tau _{g}+\kappa _{n}\tau _{g}-\kappa _{g}\tau \right] \right)
^{\prime }\right\vert _{s=l}-\left. 2\left( \frac{\kappa _{n}\tau }{\kappa
^{2}}\right) ^{\prime \prime }\right\vert _{s=l}\right\} \\
&&+\mu ^{\prime }\left( l\right) \left\{ 2\frac{\tau \left( l\right) }{%
\kappa ^{2}\left( l\right) }\left[ \kappa _{n}^{\prime }\left( l\right)
+3\kappa _{g}\left( l\right) \tau _{g}\left( l\right) +\kappa _{n}\left(
l\right) \tau _{g}\left( l\right) -\kappa _{g}\left( l\right) \tau \left(
l\right) \right] +\left. 2\left( \frac{\kappa _{n}\tau }{\kappa ^{2}}\right)
^{\prime }\right\vert _{s=l}\right\} \\
&&-2\mu ^{\prime \prime }\left( l\right) \frac{\kappa _{n}\left( l\right)
\tau \left( l\right) }{\kappa ^{2}\left( l\right) }.
\end{eqnarray*}

\bigskip In order to have%
\begin{equation*}
\left. \frac{dF}{dt}\right\vert _{t=0}=0
\end{equation*}%
for any choice of the function $\mu \left( s\right) $ satisfying $\left( \ref%
{2.5}\right) $ with arbitrary values $\mu \left( l\right) ,\ \mu ^{\prime
}\left( l\right) $ and $\mu ^{\prime \prime }\left( l\right) $ the given arc
must satisfy three boundary conditions

\begin{eqnarray}
&&\left[ \frac{\tau }{\kappa ^{2}}\left[ -\kappa _{n}\kappa ^{2}-2\tau
_{g}\kappa _{g}^{\prime }+5\kappa _{n}\tau _{g}^{2}+3\kappa _{g}\tau
_{g}^{\prime }-4\kappa _{n}\tau _{g}\tau \right. \right.  \label{2.19} \\
&&\left. \left. -\left( \frac{\tau }{\kappa ^{2}}\left[ \kappa _{n}^{\prime
}+3\kappa _{g}\tau _{g}+\kappa _{n}\tau _{g}-\kappa _{g}\tau \right] \right)
^{\prime }-\left( \frac{\kappa _{n}\tau }{\kappa ^{2}}\right) ^{\prime
\prime }\right] \right\vert _{s=l}=0  \notag
\end{eqnarray}

\begin{equation}
\left. \left[ \frac{\tau }{\kappa ^{2}}\left[ \kappa _{n}^{\prime }+3\kappa
_{g}\tau _{g}+\kappa _{n}\tau _{g}-\kappa _{g}\tau \right] +\left( \frac{%
\kappa _{n}\tau }{\kappa ^{2}}\right) ^{\prime }\right] \right\vert _{s=l}=0,
\label{2.20}
\end{equation}%
\begin{equation}
\kappa _{n}\left( l\right) \tau \left( l\right) =0  \label{2.21}
\end{equation}%
and the differential equation%
\begin{eqnarray*}
&&\kappa _{g}\tau ^{2}\left( l\right) +2\frac{\tau }{\kappa ^{2}}\left[
\kappa _{g}\kappa ^{2}\left( \tau _{g}+\tau \right) +2\kappa _{g}\kappa
_{n}\kappa _{g}^{\prime }-2\kappa _{g}\kappa _{n}^{2}\tau _{g}-\kappa
_{g}^{\prime }\tau _{g}^{\prime }+4\kappa _{n}\tau _{g}\tau _{g}^{\prime
}\right.  \\
&&\left. -2\kappa _{g}^{2}\kappa _{n}^{\prime }-2\kappa _{g}^{3}\tau
_{g}-\tau _{g}^{2}\kappa _{n}^{\prime }-2\kappa _{g}\tau _{g}^{3}+-\kappa
_{g}\tau _{g}^{\prime \prime }-2\kappa _{n}\tau _{g}^{\prime }\tau +2\kappa
_{g}\tau \left( \kappa _{g}^{2}+\kappa _{n}^{2}+\tau _{g}^{2}\right) \right] 
\\
&&-2\left[ \frac{\tau }{\kappa ^{2}}\left[ -\kappa _{n}\kappa ^{2}+5\kappa
_{n}\tau _{g}^{2}-2\tau _{g}\kappa _{g}^{\prime }+3\kappa _{g}\tau
_{g}^{\prime }-4\kappa _{n}\tau _{g}\tau \right] \right] ^{\prime }
\end{eqnarray*}%
\begin{equation}
+2\left[ \frac{\tau }{\kappa ^{2}}\left( \kappa _{n}^{\prime }+\kappa
_{n}\tau _{g}+3\kappa _{g}\tau _{g}-\kappa _{g}\tau \right) \right] ^{\prime
\prime }+2\left( \frac{\kappa _{n}\tau }{\kappa ^{2}}\right) ^{\prime \prime
\prime }=0.\ \ \ \ \ \ \ \ \ \ \ \ \ \   \label{2.22}
\end{equation}%
$\ \ \ \ \ \ \ \ $

\bigskip Observe that, any planar curve, namely a curve with zero torsion,
satisfies the above differential equation and the boundary conditions and it
is  a relaxed elastic line of second kind. Thus, we have the following
theorem:

\begin{theorem}
The intrinsic equations for a relaxed elastic line of second kind on a
connected oriented surface in Euclidean 3-space are given by the
differential equation $\left( \ref{2.22}\right) $ with the boundary
conditions $\left( \ref{2.19}\right) -\left( \ref{2.21}\right) $ at the free
end, where $\kappa _{g},\ \kappa _{n}\ $and $\tau _{g}$ are the geodesic
curvature, the normal curvature and the geodesic torsion as functions of the
arc length along the curve.

\begin{corollary}
A geodesic on an oriented surface is a relaxed elastic line of second kind
if it satisfies the differential equation%
\begin{eqnarray*}
&&\frac{\tau _{g}}{\kappa _{n}^{2}}\left( 4\kappa _{n}\tau _{g}^{\prime
}-\tau _{g}\kappa _{n}^{\prime }-2\kappa _{n}\tau _{g}^{\prime }\right)
-\left( \frac{\tau _{g}}{\kappa _{n}}\left( \tau _{g}^{2}-\kappa
_{n}^{2}\right) \right) ^{\prime } \\
&&+\left( \frac{\tau _{g}}{\kappa _{n}}\left( \kappa _{n}^{\prime }+\kappa
_{n}\tau _{g}\right) \right) ^{\prime \prime }+\left( \frac{\tau _{g}}{%
\kappa _{n}}\right) ^{\prime \prime \prime }=0,
\end{eqnarray*}%
and the boundary conditions%
\begin{equation*}
\tau _{g}\left( l\right) =\tau _{g}^{\prime }\left( l\right) =\tau
_{g}^{\prime \prime }\left( l\right) =0.
\end{equation*}
\end{corollary}

\begin{corollary}
A line of curvature on an oriented surface is a relaxed elastic line of
second kind if it satisfies the differential equation%
\begin{equation*}
\kappa _{g}\tau ^{2}\left( l\right) +2\frac{\kappa _{g}\tau }{\kappa ^{2}}%
\left( 3\kappa ^{2}\tau +2\kappa _{n}\kappa _{g}^{\prime }-2\kappa
_{g}\kappa _{n}^{\prime }\right) +\left( \kappa _{n}^{\prime }\tau +\kappa
_{n}\tau ^{\prime }\right) 
\end{equation*}%
\begin{equation*}
+2\left( \frac{\tau }{\kappa ^{2}}\left( \kappa _{n}^{\prime }-\kappa
_{g}\tau \right) \right) ^{\prime \prime }+2\left( \frac{\kappa _{n}\tau }{%
\kappa ^{2}}\right) ^{\prime \prime \prime }=0\ \ \ \ \ \ \ \ \ \ \ \ \ \ \ 
\end{equation*}%
and the boundary conditions%
\begin{equation*}
\left. \left( \frac{\tau }{\kappa ^{2}}\left( \kappa _{n}^{\prime }-\kappa
_{g}\tau \right) \right) ^{\prime }\right\vert _{s=l}+\left. \left( \frac{%
\kappa _{n}\tau }{\kappa ^{2}}\right) ^{\prime \prime }\right\vert _{s=l}=0,
\end{equation*}%
\begin{equation*}
\left. \left( \frac{\tau }{\kappa ^{2}}\left( \kappa _{n}^{\prime }-\kappa
_{g}\tau \right) \right) \right\vert _{s=l}+\left. \left( \frac{\kappa
_{n}\tau }{\kappa ^{2}}\right) ^{\prime }\right\vert _{s=l}=0,
\end{equation*}%
\begin{equation*}
\kappa _{n}\left( l\right) \tau \left( l\right) =0.
\end{equation*}
\end{corollary}

\begin{corollary}
An asymptotic curve on an oriented surface is a relaxed elastic line of
second kind if it satisfies the differential equation%
\begin{equation*}
\frac{\tau _{g}}{\kappa _{g}^{2}}\left( 2\kappa _{g}^{3}\tau _{g}-\kappa
_{g}^{\prime }\tau _{g}^{\prime }+\kappa _{g}\tau _{g}^{\prime \prime
}\right) -\left( \frac{\tau _{g}}{\kappa _{g}^{2}}\left( 3\kappa _{g}\tau
_{g}^{\prime }-2\tau _{g}\kappa _{g}^{\prime }\right) \right) ^{\prime
}+2\left( \frac{\tau _{g}^{2}}{\kappa _{g}}\right) ^{\prime \prime }=0
\end{equation*}%
and the boundary condition%
\begin{equation*}
\tau _{g}\left( l\right) =0.
\end{equation*}
\end{corollary}
\end{theorem}

\end{document}